\documentclass[12pt,a4paper]{article}
\usepackage{amssymb}
\usepackage[utf8]{inputenc}
\usepackage[english]{babel}
\usepackage{color}
\usepackage{amsmath}
\usepackage{enumerate}

\newcommand{\Z}{\mathbb{Z}}
\newcommand{\N}{\mathfrak{N}}
\newcommand{\NT}{\tilde\N}
\newcommand{\R}{\mathbb{R}}
\newcommand{\Aut}{\operatorname{Aut}}

\def\restrict#1#2{\mathchoice
{\setbox1\hbox{${\displaystyle #1}_{\scriptstyle #2}$}
\restrictaux{#1}{#2}}
{\setbox1\hbox{${\textstyle #1}_{\scriptstyle #2}$}
\restrictaux{#1}{#2}}
{\setbox1\hbox{${\scriptstyle #1}_{\scriptscriptstyle #2}$}
\restrictaux{#1}{#2}}
{\setbox1\hbox{${\scriptscriptstyle #1}_{\scriptscriptstyle #2}$}
\restrictaux{#1}{#2}}}
\def\restrictaux#1#2{{#1\,\smash{\vrule height .8\ht1 depth .85\dp1}}_{\,#2}}

\newtheorem{Thm}{Theorem}

\newtheorem{Lem}{Lemma}
\newtheorem{Cor}{Corollary}
\newtheorem{Rem}{Remark}

\title{$L^{p}$ compression of some HNN extensions}
\author{Pierre-Nicolas Jolissaint and Thibault Pillon $\;$ \footnote{Both authors supported by Swiss SNF -grant 20-137696.}}

\begin{document}

\maketitle

\begin{abstract}
In \cite{GJ}, the authors introduce a framework to prove that a large class of HNN extensions have the Haagerup property, the main motivation being Baumslag-Solitar groups. Using this framework and new tools on locally compact groups developped in \cite{CCMT}, we are able to obtain quantitative results on embeddings into Lebesgue spaces for a large class of HNN extensions. 
\end{abstract}

\section{Introduction}

In \cite{GK}, Guentner and Kaminker introduced the notion of compression exponent between metric spaces. Roughly speaking, it gives a way of quantifying how well a metric space coarsely embeds into another. We recall here the definitions. Let $X$ and $Y$ be metric spaces. A map $f : X \rightarrow Y$ is said to be \textit{large-scale Lipschitz} if we can find some constants $A,B \geq 0$ so that the following inequality holds for every $x,y\in X$:
$$
d(f(x),f(y)) \leq A d(x,y)+B.
$$
In this case, we define $R(f)$ to be the supremum over the $\alpha\in[0,1]$ such that there exist some constants $C,D\geq 0$ so that
$d(f(x),f(y))\geq Cd(x,y)^{\alpha}-D$, for every $x,y\in X$. Then, the \textit{compression exponent}, or  simply \textit{compression} of $X$ with target $Y$ is defined as $\alpha^*_{Y}(X):=\sup_{f} R(f)$, where the supremum is taken over all Large-scale Lipschitz maps $f$ from $X$ into $Y$. When $Y=L^p(\Omega)$ with $\Omega$ a standard Borel space, we set $\alpha^*_{L^p}=\alpha^*_p$.  If $G$ is a compactly generated locally compact group (e.g. a finitely generated group), we view it as a metric space with the word metric. The study of groups seen as metric spaces has offered some striking results. Let us recall one. If $G$ is a finitely generated group and if $\alpha^{\ast}(G)>0$ for some $p\in (1,+\infty)$, then $G$ satisfies the Novikov conjecture (see \cite{KY}).

In \cite{GJ} the authors introduced the notion of an $\N$-BS group. These are groups arising as HNN extensions satisfying properties similar to Baumslag-Solitar groups. In order to prove the Haagerup property for such groups, they developed a framework that we shall heavily rely on and that we now recall.

Let $\N$ be a locally compact compactly generated group and let $G$ be a closed subgroup of $\N$. Let $i_1,i_2 : H\rightarrow G$ be two inclusions of a group $H$ onto open subgroups of finite index, and assume $i_1$ and $i_2$ are conjugated by an automorphism $\varphi$ of $\N$.  The $\N$-BS group $\Gamma$ is then the HNN extension  $\operatorname{HNN}(G,H,i_1,i_2)$ whose presentation is given by $\langle S,t \vert R, ti_{1}(h)t^{-1}=i_{2}(h) \  \forall h\in H \rangle$, where $G =\langle S \vert R \rangle$. 

\begin{Thm}\label{compression}
Let $\N$ be a connected Lie group and $G$ a closed cocompact subgroup of $\N$ and let $\Gamma$ be an HNN extension as above. Then, for all $p>1$, $\alpha_{p}^{*}(\Gamma)=1$.
\end{Thm}

The strategy to prove Theorem \ref{compression} is to construct a metric space $M$ on which $\Gamma$ acts continuously, properly, cocompactly and by isometries, so that, $M$ and $\Gamma$ are quasi-isometric by Svarc Lemma: this is done in section \ref{Millefeuille}, where we also give a quantitative comparison between two natural metrics on $M$.
In the last section, we prove Theorem \ref{compression} and treat some concrete examples.

\section{The space $M$}\label{Millefeuille}

Let $\NT=\N\rtimes\Z$, where $\Z$ acts on $\N$ by iterations of $\varphi$ and let $j_\N : \Gamma\rightarrow \NT$ be the homomorphism defined by $g \mapsto (g,0)$ for $g\in G$ and $t \mapsto (1,1)$.
 Then consider $T$, the Bass-Serre tree associated with the HNN extension $\Gamma$ and denote by $j_T:\Gamma\rightarrow \operatorname{Aut}(T)$ the monomorphism induced by the action of $\Gamma$ on $T$.
For later use, recall the following result from \cite{GJ}.

\begin{Thm}\label{GJ}
The homomorphism $j:\Gamma\rightarrow\NT\times \Aut(T),\ g\mapsto (j_\N(g),j_T(g))$ is injective and has closed image. In particular, it is proper.
\hfill $\square$
\end{Thm}

Following Proposition 2.1 in \cite{CCMT}, we will define a metric space $Y$ on which $\NT$ acts continuously, properly, cocompactly and by isometries. Endow $\NT$ with a left-invariant Riemannian metric. For each coset $L_{i}=\N\times\{i\}$ of $\N$ in $\NT$ we consider a strip $L_{i}\times[0,1]$, equipped with the product Riemannian metric, and attach it to $\NT$ by identifying $(l,0)$ to $l$ and $(l,1)$ to $l\cdot (1,1)$. Denote by $Y$ the space obtained in this way. $Y$ has a natural shortest-path metric induced by the riemannian metric on each of the strips. Furthermore, Y is naturally homeomorphic (but not necessarily isometric!) to $\N\times\R$. Using this obvious parametrization, $\NT$ acts on $Y$ by $(n,k)\cdot (y,s) = (n\varphi^{k}(y), k+s)$, for $(n,k)\in \NT$ and $(y,s)\in Y$. As in Proposition $2.1$ in \cite{CCMT}, $Y$ is a locally compact, geodesic metric space on which $\NT$ acts continuously, properly, cocompactly and by isometries. We denote by $b$ the projection map $(y,s)\mapsto s$.


Let us recall briefly the construction of the Bass-Serre tree $T$ of $\Gamma$. It is an oriented graph whose vertices are the left-cosets $\Gamma / G$ and the edges correspond to the left cosets $\Gamma / i_{1}(H)$. The edge $\gamma / i_{1}(H)$ is directed from $\gamma t^{-1}G$ to $\gamma G$. As the $i_{k}(H)$ are of finite index in $G$, $T$ is locally finite. Then, by construction, $\Gamma$ acts naturally on $T$ by left multiplication.\\
Now, let $p : \Gamma \rightarrow \Z$ be the homomorphism defined on the generators by $p(t)=1$ and $p(g)=0$, for every $g\in G$. Since the vertices of $T$ correspond to the left cosets of $G$ in $\Gamma$, we can define a map $c$ on the vertices of $T$ by $c(xG)=p(x)$ and extend it to the metric tree $T$ by affine interpolation. This allows us to define the fibre product $M$\,:
$$
M=\{(x,y)\in T\times Y \,: \,c(x)=b(y)\}.
$$
The subspace $M$ is $\Gamma$-invariant for the diagonal action of $\Gamma$ on $T\times Y$. 
Indeed, for all $x\in T$, $c(g\cdot x)=c(x)$ if $g\in G$ and $c(t\cdot x)=c(x)+1$. In a similar fashion, for all $y\in Y$, $b(g\cdot y)=b(y)$ if $g\in G$ and $b(t\cdot y)=b(y)+1$. Hence, if $c(x)=b(y)$, it implies that $c(\gamma\cdot x)=b(\gamma\cdot y)$ for any $\gamma\in\Gamma$. Other similar fibre products have already been considered, namely horocyclic products and millefeuille spaces. Those spaces are defined using so-called Busemann functions (see Section 7 in \cite{CCMT}). It is worth noting that, in general, our functions $c$ or $b$ are not Busemann functions. 

We endow $T\times Y$ with the product metric, namely, $d((x,y),(x',y'))=d_{T}(x,x')+d_{Y}(y,y')$.

\begin{Lem}\label{MilfPath}
$M$ is path-connected. Furthermore, denoting by $d_{M}$ the shortest-path metric induced by $d$ on $M$, the metrics $d$ and $d_{M}$ are bilipschitz equivalent.
\end{Lem}

\noindent {\bf Proof : } First, observe that, for any point $y=(n,s)\in Y$, the path
$$
\alpha_{y} : \R \rightarrow Y : \ u\mapsto (n,u+s)
$$
is a geodesic such that $\alpha_{y}(0)=y$ and $b(\alpha_y(u))=b(y)+u$, $\forall u\in \R$. Similarly, for any point $x\in T$ one can chose a geodesic path $\beta_{x} :\R \rightarrow T$ such that $\beta_{x}(0)=x$ and $c(\beta_{x}(u))=c(x)+u$. Let $(x_{0},y_{0}),(x_{1},y_{1})\in M$. We will build a path linking those points in two steps. For the first one, let $\sigma : [0,d_{T}(x_{0},x_{1})] \rightarrow T$ be the geodesic from $x_{0}$ to $x_{1}$. Let $\theta_{1}$ be the path defined by 
$$
\theta_{1}(u)=(\sigma(u),\alpha_{y_{0}}(c(\sigma(u))-b(y_{0})).
$$ 
The left component links $x_{0}$ to $x_{1}$, while the right component starts from $y_{0}$ and ends at a certain point $y_{2}$. Moreover, the path $\theta_{1}$ is contained in $M$. Indeed, for all $u\in [0,d_{T}(x_{0},x_{1})]$, we have:
\begin{eqnarray*}
b(\alpha_{y_{0}}\left(c(\sigma(u))-b(y_{0})\right)) 
&=& 
b(y_{0})+c(\sigma(u))-b(y_{0})\\
&=&
c(\sigma(u)).
\end{eqnarray*}
So, $\theta_{1}$ connects $(x_{0},y_{0})$ to a point $(x_{1},y_{2})\in M$ satisfying $b(y_{2})=c(x_{1})=b(y_{1})$. For the second step, we will find a path in $M$ between $(x_{1},y_{2})$ and $(x_{1},y_{1})$. In a similar way, let $\tilde{\sigma} : [0,d_{Y}(y_{2},y_{1})] \rightarrow Y$ be a geodesic path linking $y_{2}$ to $y_{1}$ in $Y$. Then, it is easy to check that the path 
$$
\theta_{2} : [0,d_{Y}(y_{2},y_{1})] \rightarrow M : \theta_{2}(u)=(\beta_{x_{1}}(b(\tilde{\sigma}(u))-c(x_{1})),\tilde{\sigma}(u))
$$ 
does the job. This shows that $M$ is path-connected. Now, the inequality $d\leq d_{M}$ being immediate, we need to analyze the length of the path we just considered in order to finish the proof. Denoting by $L(\theta_{j})$ the length of the path $\theta_{j}$, we get the following estimates:
$$
L(\theta_{1})\le 2d_{T}(x_{0},x_{1})
$$
and
$$
L(\theta_{2})\le 2d_{Y}(y_{2},y_{1})\le 2(d_{Y}(y_{0},y_{1})+d_{Y}(y_{1},y_{2}))
$$
By construction, $d_{Y}(y_{0},y_{2})\leq d_{T}(x_{0},x_{1})$. We can conclude:
\begin{eqnarray*}
d_M((x_{0},y_{0}),(x_{1},y_{1}))
&\le & 
L(\theta_{1})+L(\theta_{2}) \\
&\le & 
2d_{T}(x_{0},x_{1})+2d_{Y}(y_{0},y_{1})+2d_{Y}(y_{1},y_{2}) \\
&\le &
4d_{T}(x_{0},x_{1})+2d_{Y}(y_{0},y_{1}) \\
&\le & 
4\cdot d((x_{0},x_{1}),(y_{0},y_{1})).
\end{eqnarray*}
\hfill $\square$

\section{Proof of Theorem \ref{compression} and Applications}

In order to apply Svarc Lemma, we prove that the action of $\Gamma$ is proper and cocompact.

\begin{Lem}\label{proper}
The $\Gamma$-action on $T\times Y$ is proper. That is, for all $(x,y)\in T\times Y$, there exists $r>0$ so that $\{\gamma \in \Gamma \ : \ \gamma \cdot B((x,y),r) \cap B((x,y),r) \neq \emptyset\}$ is relatively compact in $\Gamma$.
\end{Lem}
In particular, as $M$ is a closed subset of $T\times Y$, we get immediately the following Corollary.
\begin{Cor}
The $\Gamma$-action on the fibre product $M$ is proper.
\end{Cor}

\noindent {\bf Proof of Lemma \ref{proper} : } The action of $Aut(T)\times \NT$ on $T\times Y$ is proper, therefore, by Theorem \ref{GJ}, it is also the case for the action of $\Gamma$. As $M$ is closed and $\Gamma$-invariant, we can conclude. 
\hfill $\square$

\begin{Lem}
The action of $\Gamma$ on $M$ is cocompact.
\end{Lem}

\noindent{\bf Proof : } It is enough to see that, for any sequence $(x_{k},y_{k})_{k}\subset M$, we can find a sequence $(\gamma_{k})_{k}\subset\Gamma$ so that the sequence $(\gamma_{k}\cdot(x_{k},y_{k}))_{k}$ converges. Since $\Gamma$ acts transitively on the edges of $T$, we can assume that the sequence $(x_{k})_{k}$ belongs to the edge $[G,tG]$. This implies that $0 \leq c(x_{k})=b(y_{k})\leq 1$, for all but possibly finitely many $k$, so that the sequence $(y_{k})_{k}$  is contained in the strip of $Y$ corresponding to the coset $\N$ in $\NT/\N$. Using the fact that the action of $G$ on $\N$ is cocompact, we can multiply by elements of $G$ in such a way that the sequence $(y_{k})_{k}$ converges. But since $G$ stabilizes the vertex $G$ in $T$, this process maintains the sequence $(x_{k})_{k}$ inside the edges adjacent to $G$. Since there are only finitely many on these, the sequence $(x_k,y_k)_{k}$ converges up to extracting a subsequence. This concludes the proof.
\hfill $\square$ \\

We are now able to prove Theorem \ref{compression}.

\noindent{\bf Proof of Theorem \ref{compression} : } Firstly, we show that $\alpha^{\ast}_{p}(\Gamma)\geq \alpha^{\ast}_{p}(\NT)$. Indeed, by Svarc Lemma, $\Gamma$ is quasi-isometric to $(M,d_{M})$, which is quasi-isometric to $(M,d)$ by Lemma \ref{MilfPath}. Moreover, $Y$ is quasi-isometric to $\NT$. Hence, $\alpha_p^*(\Gamma)=\alpha_p^*(M)\ge \alpha^*_p(T\times Y)$ and $\alpha^*_p(Y)=\alpha^*_p(\NT)$. Then, the lower bound follows from the propositions:
\begin{itemize}
	\item For a tree $T$, $\alpha^{*}_{p}(T)=1$, for all $p>1$. (See Theorem 2.6 in \cite{BS})
	\item For two metric spaces $X$ and $X'$, the compression of $X\times X'$ is the minimum of the compressions of the factors. (See \cite{GK})
\end{itemize}
Finally, we conclude the proof by noting that $\alpha_{p}^{\ast}(\NT)=1$, which follows from the following propositions:
\begin{itemize}
	\item Any semi-direct product of a connected Lie group with $\Z$ is quasi-isometric to a connected Lie group. (By an unpublished result of Y. Cornulier)
	\item Let $\mathcal{K}$ be a connected Lie group. Then, $\alpha_{p}^{\ast}(\mathcal{K})=1$, for all $p>1$. (See \cite{T}) 

\end{itemize}

\hfill $\square$ \\

We remark that, if $\N$ is a soluble connected Lie group, then Cornulier's result is a simple consequence of a lemma of Mostow. Here is a short proof that we owe to Alain Valette. In this case, $\NT$ is soluble and Noetherian (i.e. every closed subgroup is compactly generated). A lemma of Mostow (see Lemma 5.2 in \cite{M}) asserts that there exist a compact normal subgroup $K$ of $\NT$ and a soluble almost connected Lie group $\mathcal{M}$ such that the quotient $\NT /K$ is isomorphic to $\mathcal{L}$, where $\mathcal{L}$ is a cocompact, closed subgroup of $\mathcal{M}$. Then, the connected component of unity $\mathcal{M}^{0}$ is quasi-isometric to $\NT$. Indeed, on one hand, by compactness, $\NT$ is quasi-isometric to $\NT/K$ and by Mostow we deduce that $\NT$ is quasi-isometric to $\mathcal{L}$. On the other hand, by cocompactness, $\mathcal{L}$ is quasi-isometric to $\mathcal{M}$ and, since $\mathcal{M}$ has only finitely many connected components, it is quasi-isometric to $\mathcal{M}^{0}$.\\

In particular, Theorem \ref{compression} allows us to cover all the examples appearing in \cite{GJ}.

\begin{Cor} The following groups have compression $1$.

\begin{enumerate}
	\item The Baumslag-Solitar groups $BS_{q}^{p}=\langle x,t \mid x^{p}=tx^{q}t^{-1} \rangle= \\ HNN(\Z,\Z,p\cdot, q\cdot)$, with parameters $p,q\in\Z_{+}$. (For a different proof, see also \cite{CV}.)
	\item Torsion free, finitely presented abelian-by-cyclic groups.
	\item Let $\N$ be a homogeneous nilpotent Lie group. So, it admits a dilating automorphism $\varphi$. Suppose that $\N$ contains a discrete, cocompact subgroup $G$ which is invariant by $\varphi$. Then, for any finite index subgroup $H$ in $G$, the extension $HNN(G,H,i_{1},\restrict{\varphi}{H})$, where $i_{1}$ is the canonical injection, has compression $1$.
\end{enumerate}
\end{Cor}

\noindent{\bf Proof : } The Bausmslag-Solitar groups $BS_{q}^{p}$ can be seen as HNN extension of $\Z$ with itself, considering the inclusions $i_{1}, i_{2} : \Z \rightarrow \Z $ defined by $i_{1}(n)=pn$ and $i_{2}(n)=qn$. Then, apply Theorem \ref{compression} with the automorphism of $\N=\R$ given by $\varphi(x)=\frac{p}{q}x$. \\
For the second class of examples, it is known (see for instance \cite{FM}) that these groups are HNN extensions of $\Z^{n}$ with itself with respect to $i_{1},i_{2} : \Z^{n} \rightarrow \Z^{n}$, where $i_{1}$ is the identity and $i_{2}\in GL_{n}(\Z)$. Again, apply Theorem \ref{compression} with the automorphism of $\N=\R^{n}$ given by $i_{2}^{-1}$.\\
The case of the last class of examples is clear by construction.

\hfill $\square$

\begin{Rem} 
It is important to note that the assumption about finite presentation is necessary to treat the second class of examples. Indeed, the wreath product $\Z \wr \Z$ is torsion free, abelian-by-cyclic and finitely generated. However, it is computed in \cite{ANP} that  $\alpha_{2}^{*}(\Z \wr \Z)=\frac{2}{3}$.
\end{Rem}

\begin{Rem}
In the case where $\N$ is a generic compactly generated locally compact (not necessarily a connected Lie group), it is also possible to find a space $Y$ admitting a geometric action of $\NT$, by Proposition 2.1 in \cite{CCMT}. However, it is not clear how to endow $Y$ with a natural fibration $b$ compatible with the semi-direct product structure on $\NT$ in order to generalize Theorem \ref{compression}.
\end{Rem}










\bibliography{biblio}
\bibliographystyle{amsalpha}

\bigskip

Authors addresses:
\medskip

\noindent
Institut de Math\'ematiques - Unimail\\
11 Rue Emile Argand\\
CH-2000 Neuch\^atel\\
Switzerland

\begin{verbatim}pierre-nicolas.jolissaint@unine.ch; thibault.pillon@unine.ch
\end{verbatim}

\end{document}